\documentclass[11pt]{article}
\usepackage{amsmath,amssymb}
\usepackage{latexsym}
\usepackage{url}
\usepackage{graphicx}
\usepackage{graphics}
\usepackage[bottom]{draftcopy}
\newtheorem{thm}{Theorem}[section]
\newtheorem{lem}[thm]{Lemma}
\newtheorem{result}[thm]{Result}

\newtheorem{prop}[thm]{Proposition}

\newtheorem{rmk}[thm]{Remark}

\newcommand{\epf}{\mbox{${\blacksquare}$ }}

\newcommand{\beq}{\begin{eqnarray}}
\newcommand{\eeq}{\end{eqnarray}}
\newcommand{\beqn}{\begin{eqnarray*}}
\newcommand{\eeqn}{\end{eqnarray*}}

 \newcommand{\jmod}[2]{\mbox{ $\equiv #1$}\mbox{ \rm (mod $#2$)}}

\newcommand{\njmod}[2]{\mbox{ $\not \equiv #1$}\mbox{ \rm (mod $#2$)}}

\newcommand{\bpf}{\noindent {\bf Proof:} $\mbox{}$}

 \newcommand{\CB}{\mbox{$\cal B$}}
 \newcommand{\CD}{\mbox{$\cal D$}}
 \newcommand{\CE}{\mbox{$\cal E$}}
 
 \newcommand{\CI}{\mbox{$\cal I$}}

 \newcommand{\CP}{\mbox{$\cal P$}}
 
 \newcommand{\CS}{\mbox{$\cal S$}}

 \newcommand{\R}{\mathbb{R}}
 \newcommand{\GL}{{\rm GL}}
\newcommand{\AT}{\mathbb{ATLAS}}
\newcommand{\Z}{\mathbb{Z}}

\newcommand{\F}{\mathbb{F}}
\newcommand{\la}{\langle}
\newcommand{\ra}{\rangle}

\title{A projective two-weight code related to the simple group ${\rm Co}_1$ of Conway}
\author{B.~G. Rodrigues\thanks{This work is based on the research supported
by the National Research Foundation of South Africa (Grant Numbers 95725 and 106071)}\\School of Mathematics, Statistics and Computer Science\\
University of KwaZulu-Natal \\Durban 4000, South Africa }
\date{}
\begin{document}
\maketitle
\begin{abstract}

A binary $[98280, 24, 47104]_2$ projective two-weight code related to the sporadic simple group ${\rm Co}_1$ of Conway is constructed as a faithful and absolutely irreducible submodule of the permutation module induced by the primitive action of ${\rm Co}_1$ on the cosets of ${\rm Co}_2$. The dual code of this code is a uniformly packed $[98280, 98256,3]_2$  code. The geometric significance of the codewords of the code can be traced to the vectors in the Leech lattice, thus revealing that the stabilizer of any non-zero weight codeword in the code is a maximal subgroup of ${\rm Co}_1$. Similarly, the stabilizer of the codewords of minimum weight in the dual code is a maximal subgroup of  ${\rm Co}_1$. As by-product, a new strongly regular graph on 16777216 vertices and valency 98280 is constructed using the codewords of the code.

\medskip

 \noindent {{\em
Key words and phrases\/}: strongly regular graph, flag transitive symmetric
design, automorphism group, modular representation, Conway group.
\hfill}
\end{abstract}

\section{Introduction}
Given a permutation group $G$ on a finite set $\Omega$ and a field $\F$ it is often of considerable interest to know the structure of the permutation module $\F\Omega$ (that is, the vector space over $\F$ with basis $\Omega$ considered as an $\F G$ module). The $G$-invariant submodules of $\F\Omega$ can be regarded as linear codes in $\F\Omega,$ and one may therefore ask for the weight distribution, and the partitioning of the codes into $G$-orbits (see \cite{knapp&schaeffer}). To a greater extent this paper fits into a programme outlined in \cite{knap&schmid} in that we determine codes, in particular binary codes invariant under a prescribed permutation group. However, when considering large groups, the chance of determining all codes invariant under the group decreases due to the large degree of their representation and the large dimension of the submodules of a permutation module associated to a given permutation representation. Thus, one might be satisfied with a few representations, or at least the degree of the smallest faithful (and hence irreducible) representation. Representation theory has proven to be an extremely powerful tool for the exact calculations of the minimal degrees of faithful representations of finite simple groups. In particular, see \cite{jansen} for calculations related with the sporadic simple groups and their covering groups. The dimension of irreducible representations turns out to be the $p$-rank (i.e., dimension over $\F_p$) of linear codes, (see \cite{brouvane} for an account of the minimality of the $p$-ranks of codes from combinatorial structures and \cite{moorod4,moorirodriguesco2,moorirodriguesrudv} which examine binary codes of smallest possible dimensions invariant under some sporadic simple groups). In this paper, using a modular representation theoretic approach we construct from the primitive permutation representation of degree 98280 of the simple group ${\rm Co}_1$ of Conway an irreducible faithful representation of dimension 24 as a binary code, thereby giving a nice construction of a small representation of this group. In the theorem given below, we summarize our results; the specific results relating to the codes are given as propositions and lemmas in the following sections.

\begin{thm}\label{thm-1}
Let $G$ be the simple Conway group ${\rm Co}_1$ and ${\rm C}_{24}$ be a submodule of dimension $24$ obtained from the permutation module of degree $98280$. Then the following hold:
\begin{itemize}
\item[{\rm (a)}] ${\rm C}_{24}$ is the smallest non-trivial ${\rm Co}_1$-invariant irreducible
$\F_2$-module;
\item[{\rm (b)}] ${\rm C}_{24}$ is a self-orthogonal doubly-even two-weight code;
\item[{\rm (c)}] the supports of the minimum words define uniquely a self-dual, point-primitive and flag-transitive symmetric $1$-$(98280,47104,47104)$ design invariant under ${\rm Co}_1$;
\item[{\rm (d)}] ${\rm Aut}({\rm C}_{24}) \cong {\rm Co}_1$;
\item[{\rm (e)}] the non-trivial codewords of ${\rm C}_{24}$ define a strongly regular  $(16777216, 98280, 4600, 552)$  graph $\Gamma(C_{24}).$ 
\end{itemize}
\end{thm}

The paper is organized as follows: in Section~\ref{sec-meth} we outline our background and notation and  in Section~\ref{co1group} we give a brief but complete overview on the ${\rm Co}_1$ group. In Section~\ref{bincode} we describe the construction method used and give our results on the 24-dimensional binary code invariant under ${\rm Co}_1$. In the ensuing sections, namely Sections~\ref{sec-co1des}, ~\ref{stabco1} and \ref{srgofCo1} we present our results concerning with the unique flag-transitive design on 98280 points and a strongly regular graph on 16777216 vertices invariant under the Conway group ${\rm Co}_1.$

\section{Terminology}\label{sec-meth}

We assume that the reader is familiar with some basic notions and elementary facts from design and coding theory. Our notation for codes and groups will be standard, and it is as in \cite{ak:book} and $\AT$~\cite{atlas}. For the structure of groups and their maximal subgroups we follow the $\AT$~\cite{atlas} notation.  The groups $G.H,$ $G{:}H,$ and $G^{\cdot}H$ denote a general extension, a split extension and a non-split extension respectively. For a prime $p,$ the symbol $p^m$ denotes an elementary abelian group of that order. For $G$ a finite group acting on a finite set $\Omega$, the set $\F_p\Omega$, that is, the vector space over $\F_p$ with basis $\Omega$ is called  an $\F_p G$ {\bf permutation module}, if the  action of $G$ is extended linearly on $\Omega$.

An incidence structure $\mathcal{D}=( \CP,\CB,\CI)$, with point set \CP, block set \CB\ and incidence \CI\ is a $t$-$(v,k,\lambda)$ design, if $|\CP|=v$, every block $B \in \CB$ is incident with precisely $k$ points, and every $t$ distinct points are together incident with precisely $\lambda$ blocks.  The complementary design of \CD\ is obtained by replacing all blocks of \CD\ by their complements. The design \CD\ is {\bf symmetric} if it has the same number of points and blocks. An automorphism of a design \CD\ is a permutation on \CP\ which sends blocks to blocks. The set of all automorphisms of \CD\ forms its full automorphism group denoted by Aut$\CD$.

The  {\bf code $C$ of the design} \CD\ over the finite field $\F_p$ is the space spanned by the incidence vectors of the blocks over $\F_p$. The weight enumerator of $C$ is defined as $\sum_{c \in C} x^{{\rm wt}(c)}.$ The {\bf hull} of a design \CD\ with code $C$ over the field $F$ is the code obtained by taking the intersection of $C$ and its dual. A linear $[n,k]$ code is called {\bf projective} if no two columns of a generator matrix $G$ are linearly dependent, i.e., if the columns of G are pairwise different points in a projective $(k - 1)$-dimensional space.  A {\bf two-weight} code is a code which has exactly two non-zero weights, say $w_1$ and $w_2.$ The dual of a two-weight code belongs to the important family of uniformly packed codes. A code $C$ is {\bf self-orthogonal} if $C\subseteq C^\perp$ and {\bf self-dual} if equality is attained. The all-one vector will be denoted by $\mathbf{1}$, and is the constant vector of weight the length of the code, and whose coordinate entries consist entirely of 1's. A binary code $C$ is {\bf doubly-even} if all codewords of $C$ have weight divisible by four.  Two linear codes are {\bf isomorphic} if they can be obtained from one another by permuting the coordinate positions. An {\bf automorphism} of a code is any permutation of the coordinate positions that maps codewords to codewords and will be denoted {\rm Aut}$(C).$

\section{The Conway group ${\rm Co}_1$}\label{co1group}

The {\it Leech lattice} is a certain 24-dimensional $\Z$-submodule of the 24-dimensional Euclidean space $\R^{24}$ discovered by John Leech. John Conway showed that the automorphism group of the Leech lattice is a quasisimple group. Its central factor group is the Conway group ${\rm Co}_1$. The Conway groups ${\rm Co}_2$ and ${\rm Co}_3$ are stabilizers of sublattices of the Leech lattice. We give a brief description of the construction of these groups, omitting detail. A more recent and comprehensive account is given in \cite{wilsbook}, see also \cite{conway,wilsoCo1,vectstab}.

Let $H = M_{24}$ and $(\Omega,\cal{C})$ be the Steiner system $S(24,8,5)$ for $H$. Let $V$ be the permutation module over $\F_2$ of $H$ with basis $\Omega$ and $V_{\cal{C}}$ the Golay code submodule. Let $\R^{24}$ be the permutation module over the reals for $H$ with basis $\Omega$ and let $(\phantom{.},\phantom{.})$ be the symmetric bilinear form on $\R^{24}$
for which $\Omega$ is an orthogonal basis. Then $\R^{24}$ together with $(\phantom{.},\phantom{.})$ is simply the 24-dimensional Euclidean space admitting the action of $H$, and for $\sum_{\omega}\alpha_{\omega} \omega$ and $\sum_{\omega}\beta_{\omega} \omega$ in $\R^{24}$,
\[
\left( \sum_{\omega}\alpha_{\omega} \omega, \sum_{\omega}\beta_{\omega} \omega \right) = \sum_{\omega}\alpha_\omega\beta_\omega.\]
For $v \in \R^{24}$ define $q(v) = (v,v)/16$. Thus $q$ is a positive definite quadratic form on $\R^{24}$. Given $Y \subseteq \Omega$, define $e_Y = \sum_{y \in Y}y \in \R^{24}$. For $\omega \in \Omega$ let $\lambda_\omega = e_\Omega - 4\omega$.

The Leech lattice is the set $\Lambda$ of vectors $v = \sum_{\omega}\alpha_\omega \omega \in \R^{24}$ such that:
\begin{description}
\item{($\Lambda 1$)} $\alpha_\omega \in \Z$ for all $\omega \in \Omega$.
\item{($\Lambda 2$)} $m(v) = (\sum_{\omega}\alpha_\omega)/4 \in \Z$.
\item{($\Lambda 3$)} $\alpha_\omega \jmod{m(v)}{2}$ for all $\omega \in \Omega$.
\item{($\Lambda 4$)} ${\cal{C}}(v) = \{ \omega \in \Omega \; | \; \alpha_\omega \njmod{m(v)}{4} \} \in V_{\cal{C}}$.
\end{description}
The Leech lattice $\Lambda$ is a $\Z$-submodule of $\R^{24}$. Let $\Lambda_0$ denote the set of vectors $v \in \Lambda$ such that $m(v) \jmod{0}{4}$. Then $\Lambda_0$ is a $\Z$-submodule spanned by the set $\{2e_B \; | \; B \subset \cal{C}\}$. Further, $\Lambda$ as a $\Z$-submodule is generated by $\Lambda_0$ and $\lambda_{\omega_0}$, for $\omega_0 \in \Omega$. Write $O(\R^{24})$ for the subgroup of $\GL(\R^{24})$ preserving the bilinear form $(\phantom{.},\phantom{.})$, or equivalently preserving the quadratic form $q$. Let $G$ be the subgroup of $O(\R^{24})$ acting on $\Lambda$. The group $G$ is the automorphism group of the Leech lattice. For $Y \subset \Omega$, write $\epsilon_Y$ for the element of $\GL(\R^{24})$ such that
\[ \epsilon_Y(\omega) = \left\{ \begin{array}{rl}
        -\omega & ,\mbox{if} \quad \omega \in Y, \\
        \omega & ,\mbox{if} \quad \omega \not\in Y. \end{array} \right. \]
Let $Q = \{\epsilon_Y \; | \; Y \in V_{\cal{C}}\}$. Then $K = H{\cdot}Q \leq G$. Given any positive integer $l$, write $\Lambda_l$ for the set of all vectors $v$ in $\Lambda$ with $q(v) = l$. Then $\Lambda = \cup_l \Lambda_l$. For $v = \sum_{\omega}\alpha_\omega \omega \in \Lambda$ and $i$ a non-negative integer, let
\[ S_i = \{ \omega \in \Omega \; | \; |\alpha_\omega| = i \} \;, \]
and define the shape of $v$ to be $(0^{l_0},1^{l_1},\ldots)$, where $l_i = |S_i(v)|$. Let $\Lambda_{2}^{2}$ be the set of all vectors in $\Lambda$ of shape $(2^8,0^{16})$, $\Lambda_{2}^{3}$ the vectors in $\Lambda$ of shape $(3,1^{23})$, and $\Lambda_{2}^{4}$ the vectors in $\Lambda$ of shape $(4^2,0^{22})$. Then $\Lambda_{2}^{i}$, $ 2 \leq i \leq 4$, are the orbits of $K$ on $\Lambda_2$, with $|\Lambda_{2}^{2}| = 2^7{\cdot}7594$, $|\Lambda_{2}^{3}|= 2^{12}{\cdot}24$ and $|\Lambda_{2}^{4}| = 2^2{\cdot}{\scriptsize \left( \begin{array}{c} 24 \\ 2 \end{array} \right)}$. Moreover, $|\Lambda_2| = 2^4{\cdot}3^3{\cdot}5{\cdot}7{\cdot}13$ and $K = N_G(\Lambda_{2}^{4})$.  Using this information it can be shown that $G$ acts transitively on $\Lambda_2$, $\Lambda_3$, and $\Lambda_4$. Also $K$ is a maximal subgroup of $G$ and $|G| = 2^{22}{\cdot}3^9{\cdot}5^4{\cdot}7^2{\cdot}11{\cdot}13 {\cdot}23$. Notice that $\epsilon_\Omega$ is the scalar map on $\R^{24}$ determined by $-1$, and hence is in the center of $G$. Denote by ${\rm Co}_1$ the factor group $G/\la \epsilon_\Omega \ra $. Denote by ${\rm Co}_2$ the stabilizer of a vector in $\Lambda_2$ and denote by ${\rm Co}_3$ the stabilizer of a vector in $\Lambda_3$. The groups ${\rm Co}_1$, ${\rm Co}_2$ and ${\rm Co}_3$ are the {\it Conway groups}, with $|{\rm Co}_1| = 2^{21}{\cdot}3^9{\cdot}5^4{\cdot}7^2{\cdot}11 {\cdot} 13 {\cdot}23$, $|{\rm Co}_2| = 2^{18}{\cdot}3^6{\cdot}5^3{\cdot}7{\cdot}11 {\cdot}23$ and $|{\rm Co}_3| = 2^{10}{\cdot}3^7{\cdot}5^3{\cdot}7{\cdot}11 {\cdot}23$.

In Table~\ref{fig-1} we give the primitive representations of ${\rm Co}_1$ of degree $\leq 8386560$. The first column gives the ordering of the primitive representations as given by the $\AT$~\cite{atlas} and as used in our computations;  the second gives the degrees (the number of cosets of the point stabilizer), the third the number of orbits, and the remaining columns give the size of the non-trivial orbits of the respective point stabilizers.

\begin{table}[h]
\begin{center}
\begin{tabular}{|c|c|c|c|c|c|c||c|c|c|} \hline
No. & Max. sub. & Deg. & \# & length & & & & & \\ \hline
\hline 1 & ${\rm Co}_2$ & 98280 & 4& 4600& 46575 & 47104 & && \\
\hline 2&$3^{\cdot}{\rm Suz}{:}2$ & 1545600 & 5& 5346& 22880 & 405405 & 11119682&& \\
\hline 3&$2^{11}{:}{\rm M}_{24}$ & 8292375 & 6& 3542& 48576 &1457280 &2637824 &4145152&\\
\hline 4&${\rm Co}_3$ & 8386560 & 7& 11178& 37950& 257600& 1536975 & 2608200 & 3934656  \\
\hline \hline
\end{tabular}
\caption{Maximal subgroups of ${\rm Co}_1$ of degree $\leq 8386560$}\label{fig-1}
\end{center}
\end{table}

\section{The binary $[98280, 24, 47104]_2$ code}\label{bincode}


The representations of a finite group $G$ can be constructed from permutation representations, from two representations through their Kronecker product, from representations through invariant subspaces or other methods. Iterating the constructions described above together with reductions via invariant subspaces, one obtains all irreducible representations of $G$. To find a non-trivial proper $G$-invariant subspace of $\F \Omega$ it suffices to find a vector $w \neq 0$ which lies in a proper $G$-invariant subspace $W.$ It is not difficult to see that  for $0  \neq w \in  W,$ the orbit $\{g \cdot w \,|\, g \in G\}$  spans a $G$-invariant subspace contained in $W.$ However, once $G$-invariant subspaces are determined, there remains the task of proving the irreducibility of the representations. There are a number of methods to accomplish this amongst which the Norton's irreducibility  criterion, the Meataxe (in the various flavours), spinning algorithms, to name but a few.
This approach poses an obvious computational limitation which is apparent when the permutation module has a considerably higher dimension. This is naturally true for some of the larger simple groups of sporadic type, and in particular for the first group of Conway, where the smallest primitive permutation representation is of degree 98280. It is known that the Conway group ${\rm Co}_1$ possesses a unique absolutely irreducible representation of dimension 24 over $\F_2$, (see~\cite{Co1page,liebeckPraegerSaxl,wilsbook}). However, it is not at all trivial to determine which primitive permutation representation of this group possesses this irreducible representation.

For our construction we make use of the following results:

\begin{rmk}
{\rm For $x \in {\F_q}^{n}$ and a permutation $\sigma \in S_n$ we set}
\begin{equation}\label{eqn-1}
\sigma x = (x_{\sigma^{-1}(1)}, x_{\sigma^{-1}(2)}, \ldots x_{\sigma^{-1}(n)}).
\end{equation}
\end{rmk}

\noindent Let $C$ be a linear code over $\F_q$ of length $n$ and let $G \leq {\rm Aut}(C).$ If the action of $G$ on $C$ is defined by Equation~(\ref{eqn-1}) then the code $C$ becomes an $\F_q G$-module. Note that the ambient space ${\F}^{n}_q$ is also an $\F_q G$-module with respect to the same action of $G$. We formulate the fact that $C$ is an $\F_q G$-module as the following statement. \begin{result}\label{invariance} Let $C$ be an $[n,k,d]_q$ code and let $G \leq {\rm Aut}(C).$ Then $C$ is a $k$-dimensional submodule of the ambient space ${\F}^{n}_q$, considered as an $\F_q G$-module.
\end{result}


It follows from the description provided in Section~\ref{co1group} that as an abelian group under addition, the Leech lattice $\Lambda$  is isomorphic to  $\Z^{24}$, so that $\Lambda/2\Lambda \cong {\F_2}^{24}$ is a vector space of dimension 24 of $2\cdot{\rm Co}_1$ over $\F_2$ in which a central involution acts trivially. Hence, we obtain a 24-dimensional representation of ${\rm Co}_1$ over $\F_2.$ We denote this 24-dimensional representation ${\rm C}_{24}.$ By using Result~\ref{invariance}, in Proposition~\ref{pp-1co1code} below we construct this reduction modulo 2 representation as an irreducible binary doubly-even and projective two-weight code of the permutation module of dimension 98280 over $\F_2.$  

\begin{prop}\label{pp-1co1code}
Let $G$ be the simple Conway group ${\rm Co}_1$ and ${\rm C}_{24}$ denote a submodule of dimension $24$ obtained from the permutation module of degree $98280$ over $\F_2.$  Then\\
{\rm (i)}~${\rm C}_{24}$ is a self-orthogonal doubly-even projective two-weight $[98280, 24, 47104]_2$  code with  $98280$ words of weight $47104.$\\
{\rm (ii)}~The dual code ${{\rm C}_{24}}^{\perp}$  of ${\rm C}_{24}$ is a $[98280, 98256,3]_2$  uniformly packed code with $75348000$ codewords of weight $3.$ \\
{\rm (iii)}~$\mathbf{1} \in {{\rm C}_{24}}^{\perp}$ and ${\rm C}_{24}$ is the unique submodule of its dimension on which ${\rm Co}_1$ acts absolutely irreducibly.\\
{\rm (iv)}~${\rm Aut}({\rm C}_{24}) \cong {\rm Co}_1.$
\end{prop}
\bpf (i) Identify ${\rm C}_{24}$ with the 24-dimensional vector space $\Lambda/2 \Lambda \cong {\F_2}^{24}.$ The reduction modulo 2 of the ordinary character of $2 \cdot {\rm Co}_1$ of degree 24 gives rise to a faithful $2$-modular character of ${\rm Co}_1$, see Jansen \cite[Section~4.3.17]{jansen}. Moreover, since the smallest faithful representation of $3 \cdot {\rm Suz}{:}2$ in characteristic 2 has degree 24 the minimality and hence irreducibility of ${\rm C}_{24}$ is established.  This in turn establishes the $2$-rank (dimension over $\F_2$) of ${\rm C}_{24}$. Since the $2$-rank  of ${\rm C}_{24}$  equals the dimension of the hull (i.e., $2$-rank of ${\rm C}_{24}$ equals $2$-rank of ${\rm C}_{24} \cap {{\rm C}_{24}}^{\perp}$) we deduce that ${\rm C}_{24} \subseteq {{\rm C}_{24}}^{\perp}$ and so ${\rm C}_{24}$ is self-orthogonal.\\
By \cite[Theorem A1]{vectstab}, we have that there are just three orbits of ${\rm Co}_1$ on non-zero vectors, ie, 1-dimensional spaces in $\Lambda/2 \Lambda,$ namely a single orbit on the 8386560 vectors of norm 1, obtained from the vectors of norm 6 in $\Lambda,$ and two orbits of lengths 98280 and 8292375 on the isotropic vectors in $\Lambda/2 \Lambda$ obtianed from vectors of norm 4 and 8 in $\Lambda$ respectively. These vectors are termed vectors of type 3, 2, and 4 respectively (see Section~\ref{co1group}). From this we deduce that there are only two possible weights for the non-zero codewords of ${\rm C}_{24}.$ The weight distribution of  ${\rm C}_{24}$ is given below in TABLE~2. Observe that there are exactly 98280 vectors of weight 47104, ie, minimum weight words, and these are the generating vectors of the code. Since the spanning words have weight 47104, ${\rm C}_{24}$ is doubly-even.

\begin{center}
{\footnotesize
TABLE~2: The weight distribution of ${\rm C}_{24}$\\[1ex]
\begin{tabular}{lr}
\hline \hline $i$ & $\qquad\qquad\quad\;\;\, A_i$ \\
\hline
0   & 1   \\[1ex]
47104   & 98280  \\[1ex]
49152  & 16678935 \\[1ex]
\hline \hline
\end{tabular}}
\end{center}


(ii)~Since ${\rm C}_{24} \subseteq {{\rm C}_{24}}^{\perp}$ then if $w \in {\rm C}_{24}$ it follows that $w \in{{\rm C}_{24}}^{\perp}$ and so $(w, w) =0.$ Write $w = w_1w_2\ldots w_{98280}.$ Then $\sum_{i=1}^{98280}w^{2}_i = 0.$ Since $w^{2}_i = w_i$ for all $w_i \in \F_2$ then $\sum_{i=1}^{98280}w_i =w_i\mathbf{1}$. Hence $\mathbf{1} \in {{\rm C}_{24}}^{\perp}$. Using MacWilliams identities and Pless' power moment identities we obtain the minimum weight 3 for ${{\rm C}_{24}}^{\perp}.$\footnote{We thank Markus Grassl for computing the true minimum weight for ${{\rm C}_{24}}^{\perp}$. The full weight distribution can be obtained from the author.}

Now, since ${\rm C}_{24}$ has only two non-zero weights, it is a two-weight code. Recall that a code is called projective if its dual distance is at least 3.
Moreover, a code with minimum distance $d=3$ and covering radius 2 is called uniformly packed if every vector which is not a codeword is at
distance 1 or 2 from a constant number of codewords~\cite[Theorem~7.37]{vanlitsbook}. It then follows that ${{\rm C}_{24}}^{\perp}$ is uniformly packed and thus optimal. The reader should recall that uniformly packed codes are optimal with respect to the Johnson bound.\\
(iii)~The irreducibility of ${\rm C}_{24}$ is given in the proof of (i) and the uniqueness of this submodule follows from \cite[Proposition~B, p.142]{liebeckPraegerSaxl}.\\
(iv)~Finally, since ${\rm Co}_1 \subseteq {\rm Aut}({\rm C}_{24})$ and ${\rm Co}_1$ is a primitive group of degree 98280, it follows that  ${\rm Aut}({\rm C}_{24})$ is a primitive group of degree 98280. Since there are exactly 98280 codewords of minimum weight 47104,we use this fact to determine the automorphism group of the code as a set of permutations that preserve the set of minimum weight codewords. Moreover, since the outer automorphism of  ${\rm Co}_1$ is the trivial group we deduce $ {\rm Aut}({\rm C}_{24}) \leq {\rm Co}_1$ and thus have the result.~\epf

\begin{rmk}\label{markusproof}~{\rm An alternative proof that the Conway group $G = {\rm Co}_1$ is the true automorphism group of the code follows by observing that  we have a transitive action of $G$ on both the code coordinates and the minimum weight codewords. Hence, when computing the potential additional automorphisms, we would find them in the stabilizer of a codeword (point).\\
Now, fix a minimum weight codeword $c_0$ of weight 47104. Let $\CS_1$ be the support of $c_0$, and $\CS_0$ be its complement of size 51176. Suppose that there exists an additional automorphism  $\gamma$ stabilizing $c_0.$ Then we must have $\gamma \in {\rm Sym}(\CS_0) \times {\rm Sym}(\CS_1).$
Observe  from  TABLE~2 that the supports of any two minimum weight words intersect in 23552 or 22528 coordinates. However, there is no codeword $c_1$ with support that is strictly contained in $\CS_0,$ the complement of the support of $c_0.$ Hence, there is no additional automorphism.}\end{rmk}

\section{A flag-transitive and point and block-primitive design}\label{sec-co1des}
Using a result of Key and Moori~\cite[Proposition~1]{keymoo}
and corrected in~\cite{keymoo2}, (see also ~\cite{kemoro1}) we take the supports of the codewords of minimum weight in  ${\rm
C}_{24}$ and orbit these under the action of  ${\rm Co}_1$ to form the blocks of a self-dual symmetric $1$-$(98280, 47104, 47104)$
design denoted ${\CD}_{24}$ on which ${\rm Co}_1$ acts primitively on points and on blocks.  For the sake of completeness
we state that result below.

\begin{result}\label{prop-1}
Let $G$ be a finite primitive permutation group acting on the set  $\Omega$ of size $n$. Let $\alpha \in \Omega$, and let $\Delta \not = \{\alpha\}$ be an orbit of the stabilizer $G_{\alpha}$ of $\alpha$. If $ \CB = \{ \Delta^g \mid g \in G\} $ and, given $\delta \in \Delta$, $ \CE = \{ \{\alpha,\delta\}^g \mid g \in G\}$, then $\CD=(\Omega,\CB)$ forms a symmetric $1$-$(n,|\Delta|,|\Delta|)$ design. Further, if $\Delta$ is a self-paired orbit of $G_{\alpha}$ then
$\Gamma=(\Omega,\CE)$ is a regular connected graph of valency $|\Delta|$, $\CD$ is self-dual, and $G$ acts as an automorphism group on each of these structures, primitive on vertices of the graph, and on points and blocks of the design.
\end{result}

In Proposition~\ref{pp-co1des} we examine the properties of the design ${\CD}_{24}$ and show that it is the unique flag-transitive and point-primitive design on 98280 points invariant under ${\rm Co}_1$.
\begin{prop}\label{pp-co1des}
The generating words of ${\rm C}_{24}$ form the blocks of the unique, self-dual, symmetric, flag transitive and point primitive $1$-$(98280, 47104, 47104)$ design ${\CD}_{24}$ invariant under ${\rm Co}_1$. Moreover,  ${\rm Aut}({\CD}_{24}) \cong {\rm Co}_1.$
\end{prop}

\bpf By Result~\ref{prop-1} it suffices to take the set of minimum words and orbit the images under the action of ${\rm Co}_1$ to form the blocks of a self-dual symmetric $1$-$(98280, 47104, 47104)$ design. We denote this design by ${\CD}_{24}$. The stabilizer of a point in the permutation representation of ${\rm Co}_1$ of degree 98280 is ${\rm Co}_2$ (see Table~\ref{fig-1}). From the $\AT$~~\cite[p.154]{atlas} we have that ${\rm Co}_2$ has a unique class of maximal subgroups of index 47104, namely ${\rm M^cL}$. Given a subgroup $H$ in that class, its normalizer is twice larger in ${\rm Co}_1$, meaning that there are exactly two subgroups ${\rm Co}_2$ that contain $H$. Suppose such a group exists satisfying $H < {\rm Co}_2  < {\rm Co}_1$. Then $N_{{\rm Co}_1}(H)$ contains all the elements of ${\rm Co}_1$  that fix $H$ by conjugation. Take $x \in N_{{\rm Co}_1}(H)$. Then ${{\rm Co}_2}^x$ is either ${\rm Co}_2$ or a subgroup conjugate to it. To construct all subgroups conjugate to ${\rm Co}_2$  in ${\rm Co}_1$ that contain $H$ take a transversal of $H$ in $N_{{\rm Co}_1}(H)$. This transversal has two elements, one of which is in ${\rm Co}_2$ and the other is not, showing that there is exactly one other subgroup conjugate to ${\rm Co}_2$  that contains $H.$ So there is a unique symmetric $1$-$(98280, 47104, 47104)$ design for ${\rm Co}_1$. The definition of $\Omega$ and $\CB$ respectively are inferred from Result~\ref{prop-1}, and from this it is clear that $G \subseteq {\rm Aut}(\CD_{24}).$ Once more, Result~\ref{prop-1} and the $\AT$ show that $G$ acts primitively on both $\Omega$ and $\CB$ of degree $|\Omega| = |\CB| = 98280,$ and the stabilizer of a point $\omega \in \Omega$ has exactly four orbits in $\Omega$. Hence $G_{\omega}$ fixes setwise each of $\Omega_0$, $\Omega_1$, $\Omega_3$ and $\Omega \setminus (\Omega_0 \cup \Omega_1 \cup \Omega_3) = \Omega_2$ and these are all possible $G_{\omega}$-orbits. This shows that $\CD_{24}$ is a point primitive, symmetric $1$-design. It remains to show that $G = {\rm Aut}(\CD_{24}).$ Now $G \subseteq {\rm Aut}(\CD_{24})\subseteq S_{98280},$ so ${\rm Aut}(\CD_{24})$ is a primitive permutation group on $\Omega$ of degree 98280. Moreover, ${\rm Aut}(\CD_{24})_{\omega}$ must fix $\Omega_1$ setwise, and hence ${\rm Aut}(\CD_{24})_{\omega}$ also has orbits of lengths  1, 4600, 46575, and 47104 in $\Omega$. Accordingly \cite{praegersoicher} shows that the only primitive group of degree 98280, such that ${\rm Aut}(\CD_{24})_{\omega}$ has orbit lengths 1, 4600, 46575, and 47104 is ${\rm Co}_1$.~\epf


\section{Stabilizer in ${\rm Co}_1$ of a word $w_i$ of weight $i$}\label{stabco1}
Recall that by \cite[Theorem A1]{vectstab}, there are just three orbits of ${\rm Co}_1$ on 1-dimensional spaces in $\Lambda/2 \Lambda.$ Moreover, these orbits have lengths $98280, 8292375$ and $8386560$. In Proposition~\ref{maxesC24}, we use these facts to show how the orbits split under the action of ${\rm Co}_1$ on the nonzero codewords of ${\rm C}_{24}$ (see TABLE~3) and to determine the structure of $({\rm Co}_1)_{w_i}$ where $i$ is in $W$ with $W =\{47104, 49152\}.$ For $i \in W$ we define $W_i = \{ w_i \in {\rm C}_{24} \;|\; {\rm wt}(w_i) = i\}.$ We show in Proposition~\ref{maxesC24} that $({\rm Co}_1)_{w_i}$ is a maximal subgroup of ${\rm Co}_1,$ for all $i$. Taking the support of $w_i$  and orbiting that under ${\rm Co}_1$ we form the blocks of the $1$-$(98280, i, k_i)$ designs $\CD = \CD_{w_i},$ where $k_i = |(w_i)^{{\rm Co}_1}| \times \frac{i}{98280}$. We show that ${\rm Co}_1$ acts point primitively on $\CD.$

\begin{prop}\label{maxesC24}
Let  $i \in W$ and $w_i \in W_i.$ Then $({\rm Co}_1)_{w_i}$ is a maximal subgroup of ${\rm Co}_1.$ Furthermore ${\rm Co}_1$ is primitive on $\CD_{w_i}.\ $
\end{prop}
\bpf  {\bf Case 1.} Consider $W_{47104} = \{w_i \in W\;|\; {\rm wt}(w_i) = 47104\}.$ Since $W_{47104}$ is invariant under the action of ${\rm Aut}({\rm C}_{24})$ for all $w_i\in W_{47104},$ Table 1 and Table 2 imply that ${w_i}^{{\rm Co}_1} =W_{47104}.$ Therefore $W_{47104}$ forms an orbit under the action of ${\rm Co}_1$ and thus ${\rm Co}_1$ is transitive on $W_{47104}.$  Now let $t = w_{(47104)}.$ Then $({\rm Co}_1)_t$ is a subgroup of order $2^{18}{\cdot}3^6{\cdot}5^3{\cdot}7{\cdot}11 {\cdot}23$  and thus maximal in ${\rm Co}_1$. Using Tables 1 and 2 and the orbit stabilizer theorem we deduce that $[{\rm Co}_1{:}({\rm Co}_1)_{t}] = 98280$  and so we have $({\rm Co}_1)_{t} \cong {\rm Co}_2.$ \\
{\bf Case 2.} Let $W_{49152} = \{w_i \in W \;|\; {\rm wt}(w_i) = 49152\}.$ It can be deduced from \cite[Theorem A1]{vectstab} that under the action of ${\rm Co}_1$ the set $W_{49152}$ splits into two orbits of lengths 8386560 and 8292375, say $W_{(49152)_1}$ and $W_{(49152)_2}$. Let $u = w_{(49152)_1} \in W_{(49152)_1}$ and $\overline{u} = w_{(49152)_2} \in W_{(49152)_2.}$ Then $({\rm Co}_1)_u$ is a subgroup of order $501397585920$  and thus maximal in ${\rm Co}_1$. Moreover, $({\rm Co}_1)_u \cong 2^{11}{:}{\rm M}_{24}$. (Note that there is a misprint in \cite[p. 183]{atlas} for the index $[{\rm Co}_1{:}(2^{11}{:}{\rm M}_{24})]$.) Similarly, $|({\rm Co}_1)_{\overline{u}}| = 2^{10}{\cdot}3^7{\cdot}5^3{\cdot}7{\cdot}11 {\cdot}23$, so that $({\rm Co}_1)_{\overline{u}} \cong {\rm Co}_3$.

By the transitivity of ${\rm Co}_1$ on the code coordinates, the codewords of $W_i$ form a $1$-design $\CD_{w_i}$ with $A_i$ blocks. This implies that ${\rm Co}_1$ is transitive on the blocks of $D_{w_i}$ for each $w_i$ and since $({\rm Co}_1)_{w_i}$ is a maximal subgroup of ${\rm Co}_1,$ we deduce that ${\rm Co}_1$ acts primitively on $\CD_{w_i}$ for each $i.$

\begin{center}
TABLE~3\\ Stabilizer in ${\rm Co}_1$ of a word $w_i$\\[1ex]
{\small \begin{tabular}{lrc}
\hline \hline $i$ & $\qquad\qquad\quad\;\;\, ({\rm Co}_1)_{w_i}$ & Maximality\\
\hline
0 & ${\rm Co}_1$ & No \\[1ex]
47104 & ${\rm Co}_2$ & Yes \\[1ex]
${(49152)}_{1}$ & ${\rm Co}_3$ & Yes \\[1ex]
${(49152)}_{2}$ & $2^{11}{:}{\rm M}_{24}$ & Yes \\[1ex]
\hline \hline
\end{tabular}}
\end{center}

In TABLE~4 the first column represents the codewords of weight $i$ and the second column gives the structure of the designs $\CD_{w_i}$ defined above. In the third column we list the number of blocks of $\CD_{w_i}$. We test the primitivity for the action of ${\rm Co}_1$ on $\CD_{w_i}$ in the final column.

{\small \begin{center} TABLE~4: Non-trivial point primitive $1$-designs $\CD_{w_i}$ from  ${\rm Co}_1$\\[1ex]
\begin{tabular}{llrc} \hline \hline $i$  & $ \qquad\quad\;\;\, \CD_{w_i}$ & No. of blocks & Primitivity\\ \hline
47104  &  $1$-$(98280, 47104, 47104)$ & 98280 & Yes\\[1ex]
${(49152)}_{1}$ & $1$-$(98280, 49152, 4194304)$ & 8386560 & Yes\\[1ex]
${(49152)}_{2}$  & $1$-$(98280, 49152, 4147200)$ & 8292375 & Yes\\[1ex]
\hline \hline
\end{tabular}
\end{center}~\epf}

\begin{rmk}~{\rm The geometric significance and the nature of the codewords of ${\rm C}_{24}$ and the minimum words of ${{\rm C}_{24}}^{\perp}$ can be described using the Leech lattice. This in fact is an immediate consequence of \cite[Theorem A1]{vectstab} and \cite[Theorem A2]{vectstab}.

(1). The minimum words of ${\rm C}_{24}$ are the incidence vectors of the blocks of the design ${\CD}_{24}$. They are also the minimal (type 2) vectors in the Leech lattice, see \cite[p. 183]{atlas} or \cite[p. 156]{wilsoCo1}.\\
(2). Observe (from TABLE~3) that the codewords of weight 49152 split into two classes, namely a class of codewords whose stabilizer is isomorphic to $2^{11}{:}{\rm M}_{24}$, and another with stabilizer of a codeword isomorphic to ${\rm Co}_{3}$. The class of codewords with stabilizer isomorphic to $2^{11}{:}{\rm M}_{24}$ consists of the type 4 base (or ${A_{1}^{24}}$-hole) vectors, while those vectors with stabilizer ${\rm Co}_{3}$ are known to be type 3 vectors in the Leech lattice, see  \cite[p. 183]{atlas} or \cite[p. 156]{wilsoCo1}.\\
(3).~The codewords of minimum weight 3 in ${{\rm C}_{24}}^{\perp}$ are the $222$-$S$-lattices of the Leech lattice, and they are stabilized by a group isomorphic to ${\rm U}_6(2){:}{\rm S}_3.$ This follows readily from  \cite[Theorem A2]{vectstab} or ~\cite[p. 183]{atlas}.}
\end{rmk}

\section{A strongly regular graph on 16777216 vertices related to ${\rm Co}_{1}$}\label{srgofCo1}
Since ${\rm C}_{24}$ is a two-weight code a connection can be established with strongly regular graphs. For this, let $w_1$ and $w_2$ (where $w_1 < w_2$) be the weights of a $q$-ary two-weight code $C$ of length $n$ and dimension $k.$ To $C$ we may associate a graph $\Gamma(C)$ on $q^k$ vertices as follows. The vertices of the graph are identified with the codewords and two vertices corresponding to the codewords $x$ and $y$ are adjacent if and only if $d(x,\,y) = w_1.$ From the above we obtain a new strongly regular graph $\Gamma({\rm C}_{24})$ associated to ${\rm C}_{24}$ whose properties are given in

\begin{lem}\label{srgfromC24}
$\Gamma({\rm C}_{24})$ is a strongly regular $(16777216, 98280, 4600, 552)$ graph with spectrum \\$[98280]^1,[4072]^{98280},[-24]^{16678935}$. The complementary graph  $\overline{\Gamma({\rm C}_{24})}$ of  $\Gamma({\rm C}_{24})$ is a strongly regular   $(16777216, 16678935, 16581206, 16585256)$ graph. This graph has spectrum $[16678935]^1,[23]^{16678935},$\\$[-4073]^{98280}$.
\end{lem}

\begin{rmk}

{\rm The first part of the Lemma follows using \cite[Corollary~3.7]{cardkantor}. It would be of interest to determine the $2$-rank of $\Gamma({\rm C}_{24})$ and ${\rm Aut}\Gamma({\rm C}_{24}).$  For the $2$-rank of $\Gamma({\rm C}_{24})$ and thus the dimension of the binary code ${C_2}(\Gamma({\rm C}_{24}))$ of $\Gamma({\rm C}_{24}),$ an upper bound can be deduced readily by using the spectrum of the graph. Observe from above that the eigenvalues of an adjacency matrix $A$ of $\Gamma({\rm C}_{24})$ are $\theta_{0} = 98280, \theta_{1} = 4072,$ and $\theta_{2} = -24$ with corresponding multiplicities $f_0 = 1, f_1 = 98280$ and $f_2 = 16678935.$ Now, using results of ~\cite[Section~3]{brouvane}  we obtain an upper bound on the 2-rank of $\Gamma({\rm C}_{24})$, namely that ${\rm rank}_2(\Gamma({\rm C}_{24})) \leq {\rm min}(f_1 + 1,  f_2+1) = 98281$, since $2\,|\,\theta_1-\theta_2.$}
\end{rmk}
%

\begin{center}
TABLE~6: Neighbours (distance)\\[1ex]
\begin{tabular}{ccc}
\hline \hline Distance  &  \qquad\quad\;\;\, No. of vectors & type of vector\\
\hline
0  &  98280 & 1 \\[1ex]
1 & 1 & 0\\[1ex]
 & 4600 & -\\[1ex]
 & 47104 & 1\\[1ex]
 & 46575 & -\\[1ex]
 2 & 552 & -\\[1ex]
 & 49128 & -\\[1ex]
 & 48600 & -\\[1ex]
3 & 552 & -\\[1ex]
 & 49152 & 2, 3\\[1ex]
 & 48576 & -\\[1ex]
\hline \hline
\end{tabular}
\end{center}

%


\begin{thebibliography}{30}

\bibitem{Co1page}
Atlas of finite group representations.
\newblock  http://brauer.maths.qmul.ac.uk/Atlas/v3/spor/Co1/.


\bibitem{ak:book}
E.~F. Assmus, Jr and J.~D. Key.
\newblock {\em Designs and their Codes}.
\newblock Cambridge: Cambridge University Press, 1992.
\newblock Cambridge Tracts in Mathematics, Vol.~103 (Second printing with
  corrections, 1993).

\bibitem{magma}
W. Bosma, J. Cannon, and C. Playoust.
\newblock The Magma algebra system I: The user language.
\newblock {\em J. Symbolic Comput.}, \textbf{24} (1997), 235--265.
%

\bibitem{brouvane}
A.~E. Brouwer and C.~J. van Eijl.
\newblock On the $p$-rank of the adjacency matrices of strongly regular graphs.
\newblock {\em J. Algebraic Combin.}, \textbf{1} (1992), 329--346.

\bibitem{cardkantor}
\newblock R. Calderbank and W. M. Kantor.
\newblock The geometry of two-weight codes.
\newblock {\em Bull. London Math. Soc.}, \textbf{18} (1986), 97--122.

\bibitem{conway}
J.~H. Conway.
\newblock A group of order 8,315,553,613,086,720,000.
\newblock {\em  Bull. London Math. Soc.}, \textbf{1} (1969), 79--88.

\bibitem{atlas}
J.~H. Conway, R.~T. Curtis, S.~P. Norton, R.~A. Parker, and R.~A. Wilson.
\newblock {\em An Atlas of Finite Groups}.
\newblock Oxford: Oxford University Press, 1985.

\bibitem{curtislattice}
R.~T. Curtis.
\newblock On subgroups of .O. I. Lattice stabilizers.
\newblock {\em  J. Algebra.}, \textbf{4} (1973), 549--573.

\bibitem{ganief}



\bibitem{jansen}
C.~ Jansen.
\newblock The minimal degrees of faithful representations of the sporadic simple groups and their covering groups.
\newblock  {\em LMS J. Comput. Math.}, \textbf{8} (2005), 122--144.

\bibitem{keymoo}
J.~D. Key and J.~Moori.
\newblock Designs, codes and graphs from the {J}anko groups ${J}_1$ and ${J}_2$.
\newblock  {\em J. Combin. Math. and Combin. Comput.}, \textbf{40} (2002), 143--159.

\bibitem{kemoro1}
J.~D. Key, J.~Moori, and B.~G. Rodrigues.
\newblock On some designs and codes from primitive representations of some finite simple groups.
\newblock {\em J. Combin. Math. and Combin. Comput.}, \textbf{45} (2003), 3--19.

\bibitem{keymoo2} J.~D. Key and J.~Moori.
\newblock Correction to: ``{C}odes, designs and graphs from the {J}anko groups {$J\sb 1$} and {$J\sb 2$''
[{J}. {C}ombin. {M}ath. {C}ombin. {C}omput. {\bf 40} (2002), 143--159]},
\newblock {\em J. Combin. Math. Combin. Comput.}, \textbf{64} (2008), 153


\bibitem{knap&schmid}
W. Knapp and P. Schmid.
\newblock Codes with prescribed permutation group.
\newblock {\em J. Algebra}, \textbf{67} (1980), 415--435.

\bibitem{knapp&schaeffer}
W.~ Knapp and H.~-J. Schaeffer.
\newblock On the codes related to the Higman-Sims graph.
\newblock {\em Elect. J. of Combin.}, 22 (1) (2015), \#P1.19


\bibitem{liebeckPraegerSaxl}
M. W. Liebeck, C. E. Praeger, and J. Saxl.
\newblock The maximal factorizations of the finite simple groups and their automorphism groups.
\newblock {\em Mem. Amer. Math. Soc.}, \textbf{86} (1990), no. 432,  1--151.

\bibitem{moorod4}
J.~Moori and B.~G. Rodrigues.
\newblock A self-orthogonal doubly even code invariant under {${\rm {M}^cL}:2$}.
\newblock {\em J. Combin. Theory Ser. A}, \textbf{110} (2005), no.~1, 53--69.

\bibitem{moorirodriguesco2}
J. Moori and B. G. Rodrigues.
\newblock Some designs and  codes invariant under the simple group ${\rm Co}_2.$
\newblock {\em Journal of Algebra}, {\bf 316} (2007), 649 -- 661.



\bibitem{moorirodriguesrudv}
J. Moori and B. G. Rodrigues.
\newblock Some designs and binary codes preserved by the simple group ${\rm Ru}$ of Rudvalis.
\newblock {\em Journal of Algebra}, {\bf 372} (2012), 702-710.



\bibitem{praegersoicher}
C. E. Praeger  and  L. H. Soicher.
\newblock {\em Low rank representations and graphs for sporadic groups}.
\newblock  Cambridge: Cambridge University Press. 1997.
\newblock  Australian Mathematical Society Lecture Series, Vol.~8.



\bibitem{wilsoCo1}
R. A. Wilson.
\newblock The maximal subgroups of Conway's group ${\rm Co}_1$.
\newblock {\em J. Algebra}, \textbf{85} (1983), 144--165.
\bibitem{vectstab}
R. A. Wilson.
\newblock Vector Stabilizers and Subgroups of Leech Lattice Groups.
\newblock {\em J. Algebra}, \textbf{127} (1989), 387--408.
\bibitem{wilsbook}
R.~A. Wilson.
\newblock {\em  The finite simple groups}.
\newblock  London: Springer-Verlag London Ltd., 2009.
\newblock Graduate Texts in Mathematics, Vol.~251.

\bibitem{vanlitsbook}
J.~van Lint.
\newblock {\em Introduction to Coding Theory}.
\newblock  London: Springer-Verlag London Ltd., 1999, 3rd Edition.
\newblock Graduate Texts in Mathematics, Vol.~86.
\end{thebibliography}
\end{document}